\documentclass[11pt]{article}

\usepackage{amssymb,latexsym,amsmath}

\usepackage{graphicx}

\begin{document}

\newcommand{\bfi}{\bfseries\itshape}

\makeatletter

\@addtoreset{figure}{section}

\def\thefigure{\thesection.\@arabic\c@figure}

\def\fps@figure{h, t}

\@addtoreset{table}{bsection}

\def\thetable{\thesection.\@arabic\c@table}

\def\fps@table{h, t}

\@addtoreset{equation}{section}

\def\theequation{\thesubsection.\arabic{equation}}

\makeatother

\newtheorem{thm}{Theorem}[section]

\newtheorem{prop}[thm]{Proposition}

\newtheorem{lema}[thm]{Lemma}

\newtheorem{cor}[thm]{Corollary}

\newtheorem{defi}[thm]{Definition}

\newtheorem{convention}[thm]{Convention}

\newtheorem{rk}[thm]{Remark}

\newtheorem{exempl}{Example}[section]

\newenvironment{exemplu}{\begin{exempl}  \em}{\hfill $\square$

\end{exempl}}

\newcommand{\comment}[1]{\par\noindent{\raggedright\texttt{#1}

\par\marginpar{\textsc{Comment}}}}

\newcommand{\todo}[1]{\vspace{5 mm}\par \noindent \marginpar{\textsc{ToDo}}\framebox{\begin{minipage}[c]{0.95 \textwidth}

\tt #1 \end{minipage}}\vspace{5 mm}\par}

\newcommand{\ea}{\mbox{{\bf a}}}

\newcommand{\eu}{\mbox{{\bf u}}}

\newcommand{\ueu}{\underline{\eu}}

\newcommand{\ueo}{\overline{u}}

\newcommand{\oeu}{\overline{\eu}}

\newcommand{\ew}{\mbox{{\bf w}}}

\newcommand{\ef}{\mbox{{\bf f}}}

\newcommand{\eF}{\mbox{{\bf F}}}

\newcommand{\eC}{\mbox{{\bf C}}}

\newcommand{\en}{\mbox{{\bf n}}}

\newcommand{\eT}{\mbox{{\bf T}}}

\newcommand{\eL}{\mbox{{\bf L}}}

\newcommand{\eR}{\mbox{{\bf R}}}

\newcommand{\eV}{\mbox{{\bf V}}}

\newcommand{\eU}{\mbox{{\bf U}}}

\newcommand{\ev}{\mbox{{\bf v}}}

\newcommand{\eve}{\mbox{{\bf e}}}

\newcommand{\uev}{\underline{\ev}}

\newcommand{\eY}{\mbox{{\bf Y}}}

\newcommand{\eK}{\mbox{{\bf K}}}

\newcommand{\eP}{\mbox{{\bf P}}}

\newcommand{\eS}{\mbox{{\bf S}}}

\newcommand{\eJ}{\mbox{{\bf J}}}

\newcommand{\eB}{\mbox{{\bf B}}}

\newcommand{\eH}{\mbox{{\bf H}}}

\newcommand{\leb}{\mathcal{ L}^{n}}

\newcommand{\eI}{\mathcal{ I}}

\newcommand{\eE}{\mathcal{ E}}

\newcommand{\hen}{\mathcal{H}^{n-1}}

\newcommand{\eBV}{\mbox{{\bf BV}}}

\newcommand{\eA}{\mbox{{\bf A}}}

\newcommand{\eSBV}{\mbox{{\bf SBV}}}

\newcommand{\eBD}{\mbox{{\bf BD}}}

\newcommand{\eSBD}{\mbox{{\bf SBD}}}

\newcommand{\ecs}{\mbox{{\bf X}}}

\newcommand{\eg}{\mbox{{\bf g}}}

\newcommand{\paromega}{\partial \Omega}

\newcommand{\gau}{\Gamma_{u}}

\newcommand{\gaf}{\Gamma_{f}}

\newcommand{\sig}{{\bf \sigma}}

\newcommand{\gac}{\Gamma_{\mbox{{\bf c}}}}

\newcommand{\deu}{\dot{\eu}}

\newcommand{\dueu}{\underline{\deu}}

\newcommand{\dev}{\dot{\ev}}

\newcommand{\duev}{\underline{\dev}}

\newcommand{\weak}{\stackrel{w}{\approx}}

\newcommand{\mild}{\stackrel{m}{\approx}}

\newcommand{\strong}{\stackrel{s}{\approx}}

\newcommand{\weakdown}{\rightharpoondown}

\newcommand{\opg}{\stackrel{\mathfrak{g}}{\cdot}}

\newcommand{\opunu}{\stackrel{1}{\cdot}}
\newcommand{\opdoi}{\stackrel{2}{\cdot}}

\newcommand{\opn}{\stackrel{\mathfrak{n}}{\cdot}}
\newcommand{\opx}{\stackrel{x}{\cdot}}

\newcommand{\tr}{\ \mbox{tr}}

\newcommand{\Ad}{\ \mbox{Ad}}

\newcommand{\ad}{\ \mbox{ad}}

\title{Linear dilatation structures and inverse semigroups}

\author{Marius Buliga \\
\\
Institute of Mathematics, Romanian Academy \\
P.O. BOX 1-764, RO 014700\\
Bucure\c sti, Romania\\
{\footnotesize Marius.Buliga@imar.ro}}

\date{This version:  07.06.2007}

\maketitle

\begin{abstract} 

A dilatation structure encodes the approximate self-similarity of a 
metric space. 
A metric space $(X,d)$ which admits a strong dilatation structure (definition 
\ref{defweakstrong}) has a metric tangent space at any point $x \in X$ (theorem 
\ref{thcone}), and any such metric tangent space has an algebraic structure 
of a conical group (theorem \ref{tgene}). Particular examples of conical  
groups are Carnot groups: these are simply connected Lie 
groups  whose Lie algebra admits a positive graduation.

The dilatation structures associated to conical (or Carnot) groups are linear, 
in the sense of definition \ref{defilin}. 
Thus conical groups are the right generalization of 
normed vector spaces, from the point of view of dilatation structures.

 We prove that for dilatation structures linearity is equivalent to a 
 statement about the inverse semigroup generated by the family of dilatations 
 forming a dilatation structure on a metric space.
 
The result is new for Carnot groups and the proof seems to be new even for the
particular case of normed vector spaces.
\end{abstract}


\paragraph{MSC classes:} 	20M18; 22E20; 20F65


\section{Inverse semigroups and Menelaos theorem}

\begin{defi}
A semigroup $S$ is an inverse semigroup if for any $x \in S$ there is 
an unique element $\displaystyle x^{-1} \in S$ such that 
$\displaystyle x \, x^{-1}  x = x$ and $\displaystyle x^{-1} x \,  x^{-1}  =
x^{-1}$. 
\label{dinvsemi}
\end{defi}

An important  example of an inverse semigroup is $I(X)$, the class of all
bijective maps $\phi: \, dom \, \phi \, \rightarrow \, im \, \phi$, with 
$dom \, \phi , \, im \, \phi \subset X$. The semigroup operation is the
composition of functions in the largest domain where this makes sense. 

By the Vagner-Preston representation theorem \cite{howie} every inverse
semigroup is isomorphic to a subsemigroup of $I(X)$, for some set $X$.

\subsection{A toy example}

Let  $(\mathbb{V}, \| \cdot \|)$ be a finite dimensional, normed, real vector 
space.  By definition 
the dilatation based at $x$, of coefficient $\varepsilon>0$, is the function 
$$\delta^{x}_{\varepsilon}: \mathbb{V} \rightarrow \mathbb{V} \quad , \quad 
\delta^{x}_{\varepsilon} y = x + \varepsilon (-x+y) \quad . $$
For fixed $x$ the dilatations based at $x$ form a one parameter group which
 contracts any bounded neighbourhood of $x$ to a point, uniformly with respect 
to $x$. 

With the distance $d$ induced by the norm, the metric space $(\mathbb{V}, d)$ 
is complete and locally compact. For any $x \in \mathbb{V}$ and any 
$\varepsilon > 0$ the distance $d$ behaves well with respect to the 
dilatation $\displaystyle \delta^{x}_{\varepsilon}$ in the sense: for any 
$u, v \in \mathbb{V}$ we have 
\begin{equation}
\frac{1}{\varepsilon} \, d(\delta^{x}_{\varepsilon} u , 
\delta^{x}_{\varepsilon} v) \ = \  
\, d(u,v) \quad . 
\label{coned}
\end{equation}

Dilatations encode much more than the metric structure of the space 
$(\mathbb{V}, d)$. Indeed, we can reconstruct the algebraic structure of the 
vector space $\mathbb{V}$ from dilatations. For example let us define for any 
$\displaystyle x,u,v \in \mathbb{V}$ and $\varepsilon>0$:  
$$\Sigma_{\varepsilon}^{x}(u,v) = 
\delta_{\varepsilon^{-1}}^{x} \delta_{\varepsilon}^{\delta_{\varepsilon}^{x} u}
 (v) \quad . $$
 A simple computation shows that $\displaystyle 
 \Sigma_{\varepsilon}^{x}(u,v) =   u+ \varepsilon(-u+x) + (-x+v)$, therefore we
 can recover the addition operation in $\mathbb{V}$ by using the formula: 
\begin{equation}
\lim_{\varepsilon\rightarrow 0} \Sigma_{\varepsilon}^{x}(u,v) =  u+(-x+v)
\quad .
\label{exad}
\end{equation}
This is the  addition operation 
translated such that the neutral element is $x$. Thus, for $x=0$, we recover 
the usual addition   operation.

Affine continuous  transformations $A:\mathbb{V} \rightarrow \mathbb{V}$ admit the following 
description in terms of dilatations. 
A continuous transformation  $A:\mathbb{V} \rightarrow \mathbb{V}$ is affine if and only if for any 
$\varepsilon \in  (0,1)$, $x,y \in \mathbb{V}$ we have 
\begin{equation}
 A \,  \delta_{\varepsilon}^{x} y \ = \ \delta_{\varepsilon}^{Ax} Ay \quad . 
 \label{eq1proplin} 
 \end{equation}
Any dilatation is an affine transformation, hence for any 
 $x, y \in \mathbb{V}$ and $\varepsilon, \mu > 0$ we have 
 \begin{equation}
 \delta^{y}_{\mu} \, \delta^{x}_{\varepsilon} \ = \
 \delta^{\delta^{y}_{\mu}x}_{\varepsilon} \delta^{y}_{\mu} \quad . 
 \label{lindil}
 \end{equation}
Moreover, some compositions of dilatations are dilatations. This is precisely
stated in the next theorem, which is equivalent with the Menelaos theorem in 
euclidean geometry. 

\begin{thm}
For any  $x, y \in \mathbb{V}$ and $\varepsilon, \mu > 0$ such that $\varepsilon
\mu \not = 1$ there exists an unique
$w \in \mathbb{V}$ such that 
$$\delta^{y}_{\mu} \, \delta^{x}_{\varepsilon} \ = \ 
\delta^{w}_{\varepsilon \mu} \quad . $$
\label{teunu}
\end{thm}

For the proof see Artin \cite{artin}. A straightforward consequence of this theorem is the
following result. 

\begin{cor}
The inverse subsemigroup of $I(\mathbb{V})$ generated by  dilatations of the 
space $\mathbb{V}$ is made of all dilatations and all translations in
$\mathbb{V}$. 
\label{corunu}
\end{cor}

\paragraph{Proof.}
Indeed, by theorem \ref{teunu} a composition of two dilatations with
coefficients $\varepsilon, \mu$ with $\varepsilon \mu \not = 1$ is a dilatation.
By direct computation, if $\varepsilon \mu = 1$ then we obtain
translations. This is in fact compatible with (\ref{exad}), but is a stronger
statement, due to the fact that dilatations are affine in the sense of relation 
(\ref{lindil}).  

Moreover any translation can be expressed as a composition of two dilatations 
with coefficients $\varepsilon, \mu$ such that  $\varepsilon \mu  = 1$. Finally,
any composition between a translation and a dilatation is again a dilatation.  \quad $\square$

\subsection{Focus on dilatations}

Suppose that we take the dilatations as basic data for the toy example above. 
Namely, instead of giving to the space $\mathbb{V}$ a structure of real, normed 
vector space, we give only the distance $d$ and the dilatations $\displaystyle 
\delta^{x}_{\varepsilon}$ for all $x \in X$ and $\varepsilon > 0$. We should add
some relations which prescribe: 
\begin{enumerate}
\item[-] the behaviour of the distance 
with respect to dilatations, for example some   form of relation 
(\ref{coned}), 
\item[-] the interaction between dilatations, for example the existence of 
the limit from the left hand side of relation (\ref{exad}). 
\end{enumerate}
We denote such a collection of data by $(\mathbb{V},d, \delta)$ and call it 
a dilatation structure (see further definition \ref{defweakstrong}). 

In this paper we ask if there is any relationship between dilatations and 
inverse semigroups, generalizing relation (\ref{lindil}) and corollary
\ref{corunu}. 

Dilatation structures are far more general than our toy example. A dilatation structure on a metric space,  introduced in \cite{buligadil1}, 
is a notion in between a group and  a differential structure, expressing 
the approximate self-similarity of the metric space where it lives.

A metric space $(X,d)$ which admits a strong dilatation structure (definition 
\ref{defweakstrong}) has a metric tangent space at any point $x \in X$ (theorem 
\ref{thcone}), and any such metric tangent space has an algebraic structure 
of a conical group (theorem \ref{tgene}). Conical groups are particular examples
of contractible groups. An important class of of conical groups is formed by 
Carnot groups: these are simply connected Lie 
groups  whose Lie algebra admits a positive graduation. Carnot groups appear in many situations, in particular  in relation with
sub-riemannian geometry cf. Bella\"{\i}che   \cite{bell}, groups with 
polynomial growth cf. Gromov \cite{gromovgr}, or Margulis type rigidity  
results cf. Pansu \cite{pansu}. 

The dilatation structures associated to conical (or Carnot) groups are linear, 
in the sense of relation (\ref{lindil}), see also definition \ref{defilin}. We
actually proved in \cite{buligacont} (here theorem \ref{tdilatlin}) that a 
linear dilatation structure always comes from some
associated conical group. Thus conical groups are the right generalization of 
normed vector spaces, from the point of view of dilatation structures.

\section{Dilatation structures}

We present here an introduction into the subject of dilatation 
structures, following Buliga \cite{buligadil1}.

\subsection{Notations}

Let $\Gamma$ be  a topological separated commutative group  endowed with a continuous group morphism 
$$\nu : \Gamma \rightarrow (0,+\infty)$$ with $\displaystyle \inf \nu(\Gamma)  =  0$. Here $(0,+\infty)$ is 
taken as a group with multiplication. The neutral element of $\Gamma$ is denoted by $1$. We use the multiplicative notation for the operation in $\Gamma$. 

The morphism $\nu$ defines an invariant topological filter on $\Gamma$ (equivalently, an end). Indeed, 
this is the filter generated by the open sets $\displaystyle \nu^{-1}(0,a)$, $a>0$. From now on 
we shall name this topological filter (end) by "0" and we shall write $\varepsilon \in \Gamma \rightarrow 
0$ for $\nu(\varepsilon)\in (0,+\infty) \rightarrow 0$. 

The set $\displaystyle \Gamma_{1} = \nu^{-1}(0,1] $ is a semigroup. We note $\displaystyle 
\bar{\Gamma}_{1}= \Gamma_{1} \cup \left\{ 0\right\}$
On the set $\displaystyle 
\bar{\Gamma}= \Gamma \cup \left\{ 0\right\}$ we extend the operation on $\Gamma$ by adding the rules  
$00=0$ and $\varepsilon 0 = 0$ for any $\varepsilon \in \Gamma$. This is in agreement with the invariance 
of the end $0$ with respect to translations in $\Gamma$.

The space $(X,d)$ is a complete, locally compact metric space. For any $r>0$  
and any $x \in X$ we denote by $B(x,r)$ the open ball of center $x$ and radius 
$r$ in the metric space $X$.

On the metric space $(X,d)$ we work with the topology (and uniformity) induced 
by the distance. For any $x \in X$ we denote by $\mathcal{V}(x)$ the topological
filter of open neighbourhoods of $x$.

\subsection{Axioms of dilatation structures}

The first axiom is  a preparation for the next axioms. That is why we 
counted it as axiom 0.

\begin{enumerate}
\item[{\bf A0.}] The dilatations $$ \delta_{\varepsilon}^{x}: U(x) 
\rightarrow V_{\varepsilon}(x)$$ are defined for any 
$\displaystyle \varepsilon \in \Gamma, \nu(\varepsilon)\leq 1$. The sets 
$\displaystyle U(x), V_{\varepsilon}(x)$ are open neighbourhoods of $x$.  
All dilatations are homeomorphisms (invertible, continuous, with 
continuous inverse). 

We suppose  that there is a number  $1<A$ such that for any $x \in X$ we have 
$$\bar{B}_{d}(x,A) \subset U(x)  \ .$$
 We suppose that for all $\varepsilon \in \Gamma$, $\nu(\varepsilon) \in 
(0,1)$, we have 
$$ B_{d}(x,\nu(\varepsilon)) \subset \delta_{\varepsilon}^{x} B_{d}(x,A) 
\subset V_{\varepsilon}(x) \subset U(x) \ .$$

 There is a number $B \in (1,A]$ such that  for 
 any $\varepsilon \in \Gamma$ with $\nu(\varepsilon) \in (1,+\infty)$ the 
 associated dilatation  
$$\delta^{x}_{\varepsilon} : W_{\varepsilon}(x) \rightarrow B_{d}(x,B) \ , $$
is injective, invertible on the image. We shall suppose that 
$\displaystyle  W_{\varepsilon}(x) \in \mathcal{V}(x)$, that   
$\displaystyle V_{\varepsilon^{-1}}(x) \subset W_{\varepsilon}(x) $
and that for all $\displaystyle \varepsilon \in \Gamma_{1}$ and 
$\displaystyle u \in U(x)$ we have
$$\delta_{\varepsilon^{-1}}^{x} \ \delta^{x}_{\varepsilon} u \ = \ u \ .$$
\end{enumerate}

We have therefore  the following string of inclusions, for any $\varepsilon \in \Gamma$, $\nu(\varepsilon) \leq 1$, and any $x \in X$:
$$ B_{d}(x,\nu(\varepsilon)) \subset \delta^{x}_{\varepsilon}  B_{d}(x, A) 
\subset V_{\varepsilon}(x) \subset 
W_{\varepsilon^{-1}}(x) \subset \delta_{\varepsilon}^{x}  B_{d}(x, B) \quad . $$

A further technical condition on the sets  $\displaystyle V_{\varepsilon}(x)$ and $\displaystyle W_{\varepsilon}(x)$  will be given just before the axiom A4. (This condition will be counted as part of 
axiom A0.)

\begin{enumerate}
\item[{\bf A1.}]  We  have 
$\displaystyle  \delta^{x}_{\varepsilon} x = x $ for any point $x$. We also have $\displaystyle \delta^{x}_{1} = id$ for any $x \in X$.

Let us define the topological space
$$ dom \, \delta = \left\{ (\varepsilon, x, y) \in \Gamma \times X \times X 
\mbox{ : } \quad \mbox{ if } \nu(\varepsilon) \leq 1 \mbox{ then } y \in U(x) \,
\, , 
\right.$$ 
$$\left. \mbox{  else } y \in W_{\varepsilon}(x) \right\} $$ 
with the topology inherited from the product topology on 
$\Gamma \times X \times X$. Consider also $\displaystyle Cl(dom \, \delta)$, 
the closure of $dom \, \delta$ in $\displaystyle \bar{\Gamma} \times X \times X$ with product topology. 
The function $\displaystyle \delta : dom \, \delta \rightarrow  X$ defined by 
$\displaystyle \delta (\varepsilon,  x, y)  = \delta^{x}_{\varepsilon} y$ is continuous. Moreover, it can be continuously extended to $\displaystyle Cl(dom \, \delta)$ and we have 
$$\lim_{\varepsilon\rightarrow 0} \delta_{\varepsilon}^{x} y \, = \, x \quad . $$

\item[{\bf A2.}] For any  $x, \in K$, $\displaystyle \varepsilon, \mu \in \Gamma_{1}$ and $\displaystyle u \in 
\bar{B}_{d}(x,A)$   we have: 
$$ \delta_{\varepsilon}^{x} \delta_{\mu}^{x} u  = \delta_{\varepsilon \mu}^{x} u  \ .$$

\item[{\bf A3.}]  For any $x$ there is a  function $\displaystyle (u,v) \mapsto d^{x}(u,v)$, defined for any $u,v$ in the closed ball (in distance d) $\displaystyle 
\bar{B}(x,A)$, such that 
$$\lim_{\varepsilon \rightarrow 0} \quad \sup  \left\{  \mid \frac{1}{\varepsilon} d(\delta^{x}_{\varepsilon} u, \delta^{x}_{\varepsilon} v) \ - \ d^{x}(u,v) \mid \mbox{ :  } u,v \in \bar{B}_{d}(x,A)\right\} \ =  \ 0$$
uniformly with respect to $x$ in compact set. 

\end{enumerate}

\begin{rk}
The "distance" $d^{x}$ can be degenerated: there might exist  
$\displaystyle v,w \in U(x)$ such that $\displaystyle d^{x}(v,w) = 0$. 
\label{imprk}
\end{rk}

For  the following axiom to make sense we impose a technical condition on the co-domains $\displaystyle V_{\varepsilon}(x)$: for any compact set $K \subset X$ there are $R=R(K) > 0$ and 
$\displaystyle \varepsilon_{0}= \varepsilon(K) \in (0,1)$  such that  
for all $\displaystyle u,v \in \bar{B}_{d}(x,R)$ and all $\displaystyle \varepsilon \in \Gamma$, $\displaystyle  \nu(\varepsilon) \in (0,\varepsilon_{0})$,  we have 
$$\delta_{\varepsilon}^{x} v \in W_{\varepsilon^{-1}}( \delta^{x}_{\varepsilon}u) \ .$$

With this assumption the following notation makes sense:
$$\Delta^{x}_{\varepsilon}(u,v) = \delta_{\varepsilon^{-1}}^{\delta^{x}_{\varepsilon} u} \delta^{x}_{\varepsilon} v . $$
The next axiom can now be stated: 
\begin{enumerate}
\item[{\bf A4.}] We have the limit 
$$\lim_{\varepsilon \rightarrow 0}  \Delta^{x}_{\varepsilon}(u,v) =  \Delta^{x}(u, v)  $$
uniformly with respect to $x, u, v$ in compact set. 
\end{enumerate}

\begin{defi}
A triple $(X,d,\delta)$ which satisfies A0, A1, A2, A3, but $\displaystyle d^{x}$ is degenerate for some 
$x\in X$, is called degenerate dilatation structure. 

If the triple $(X,d,\delta)$ satisfies A0, A1, A2, A3 and 
 $\displaystyle d^{x}$ is non-degenerate for any $x\in X$, then we call it  a 
 dilatation structure. 

 If a  dilatation structure satisfies A4 then we call it strong dilatation 
 structure. 
 \label{defweakstrong}
\end{defi}
 
\section{Normed conical groups}

We shall need further the notion of normed conical group. 
Motivated by the 
case of a Lie group endowed with a Carnot-Carath\'eodory  distance induced 
by a left invariant distribution, we shall use the following definition of a
local uniform group. 

Let $G$ be a group. We introduce first the double of $G$, as the group $G^{(2)} \ = \ G \times G$ 
with operation
$$(x,u) (y,v) \ = \ (xy, y^{-1}uyv) \quad .$$
The operation on the group $G$, seen as the function
$\displaystyle op: G^{(2)} \rightarrow G$, $\displaystyle op(x,y)  =  xy$
is a group morphism. Also the inclusions:
$$i': G \rightarrow G^{(2)} \ , \ \ i'(x) \ = \ (x,e) $$
$$i": G \rightarrow G^{(2)} \ , \ \ i"(x) \ = \ (x,x^{-1}) $$
are group morphisms.

\begin{defi}
\begin{enumerate}
\item[1.]
$G$ is an uniform group if we have two uniformity structures, on $G$ and
$G\times G$,  such that $op$, $i'$, $i"$ are uniformly continuous.

\item[2.] A local action of a uniform group $G$ on a uniform  pointed space $(X, x_{0})$ is a function
$\phi \in W \in \mathcal{V}(e)  \mapsto \hat{\phi}: U_{\phi} \in \mathcal{V}(x_{0}) \rightarrow
V_{\phi}  \in \mathcal{V}(x_{0})$ such that:
\begin{enumerate}
\item[(a)] the map $(\phi, x) \mapsto \hat{\phi}(x)$ is uniformly continuous from $G \times X$ (with product uniformity)
to  $X$,
\item[(b)] for any $\phi, \psi \in G$ there is $D \in \mathcal{V}(x_{0})$
such that for any $x \in D$ $\hat{\phi \psi^{-1}}(x)$ and $\hat{\phi}(\hat{\psi}^{-1}(x))$ make sense and   $\hat{\phi \psi^{-1}}(x) = \hat{\phi}(\hat{\psi}^{-1}(x))$.
\end{enumerate}

\item[3.] Finally, a local group is an uniform space $G$ with an operation defined in a neighbourhood of $(e,e) \subset G \times G$ which satisfies the uniform group axioms locally.
\end{enumerate}
\label{dunifg}
\end{defi}

\begin{defi} A normed (local) conical group  $(G, \delta, \| \cdot \|)$ is  (local) group endowed with: (I)  a (local) action of
$\Gamma$ by morphisms $\delta_{\varepsilon}$ such that
$\displaystyle \lim_{\varepsilon \rightarrow 0} \delta_{\varepsilon} x \ = \ e$ for any
$x$ in a neighbourhood of the neutral element $e$; (II) a continuous norm  
function $\displaystyle \|\cdot \| : G \rightarrow \mathbb{R}$ which satisfies 
(locally, in a neighbourhood  of the neutral element $e$) the properties: 
 \begin{enumerate}
 \item[(a)] for any $x$ we have $\| x\| \geq 0$; if $\| x\| = 0$ then $x=e$, 
 \item[(b)] for any $x,y$ we have $\|xy\| \leq \|x\| + \|y\|$, 
 \item[(c)] for any $x$ we have $\displaystyle \| x^{-1}\| = \|x\|$, 
 \item[(d)] for any $\varepsilon \in \Gamma$, $\nu(\varepsilon) \leq 1$ and any 
 $x$  we have 
 $\displaystyle  \| \delta_{\varepsilon} x \| = \, \nu(\varepsilon) \,  \| x \|$. 
  \end{enumerate}
  \label{dnco}
  \end{defi}

Particular cases of normed conical groups are: 
\begin{enumerate}
\item[-] Carnot groups, that is simply 
 connected real Lie groups whose Lie algebra admits a positive graduation, 
 \item[-]  nilpotent p-adic groups admitting a contractive automorphism. 
 \end{enumerate}
 
 A very particular case of a 
normed conical group is described in the toy example: to any real, finite
dimensional, normed  vector space $\mathbb{V}$ we may associate the normed
conical group $(\mathbb{V},+, \delta, \| \cdot \|)$, with dilatations $\delta$ 
previously described.

 In a normed conical group $(G, \delta)$  we define dilatations based in any 
 point $x \in G$ by 
 \begin{equation}
 \delta^{x}_{\varepsilon} u = x \delta_{\varepsilon} ( x^{-1}u)  . 
 \label{dilat}
 \end{equation}
 There is also a natural left invariant distance given by
\begin{equation}
d(x,y) = \| x^{-1}y\| \quad . 
\label{dnormed}
\end{equation}
 The following result is theorem 15 \cite{buligadil1}. 

\begin{thm}
Let $(G, \delta, \| \cdot \|)$ be  a locally compact  normed group with 
dilatations. Then $(G, \delta, d)$ is 
a strong  dilatation structure, where $\delta$ are the dilatations defined by (\ref{dilat}) and the distance $d$ is induced by the norm as in (\ref{dnormed}). 
\label{tgrd}
\end{thm}

\section{Properties of dilatation structures}
\label{induced}

The following two theorems describe the most important metric and algebraic 
properties of a dilatation structure. As presented 
here these are condensed statements, available in full length as theorems 7, 8,
10 in \cite{buligadil1}.

\begin{thm}
Let $(X,d,\delta)$ be a    dilatation structure. Then the metric space $(X,d)$ 
admits a metric tangent space at $x$, for any point $x\in X$. 
More precisely we have  the following limit: 
$$\lim_{\varepsilon \rightarrow 0} \ \frac{1}{\varepsilon} \sup \left\{  \mid d(u,v) - d^{x}(u,v) \mid \mbox{ : } d(x,u) \leq \varepsilon \ , \ d(x,v) \leq \varepsilon \right\} \ = \ 0 \ .$$
\label{thcone}
\end{thm}

\begin{thm}
Let $(X,d,\delta)$ be a strong dilatation structure. Then for any $x \in X$ the triple 
 $\displaystyle (U(x), \Sigma^{x}, \delta^{x}, d^{x})$ is a normed local 
 conical group. This means: 
 \begin{enumerate}
 \item[(a)]  $\displaystyle \Sigma^{x}$ is a local group operation on $U(x)$,
 with $x$ as neutral element and $\displaystyle \, inv^{x}$ as the inverse element
 function; 
  \item[(b)] the distance $\displaystyle d^{x}$ is left invariant with respect 
  to the group operation from point (a); 
 \item[(c)] For any $\varepsilon \in \Gamma$, $\nu(\varepsilon) \leq 1$, the 
 dilatation $\displaystyle \delta^{x}_{\varepsilon}$ is an automorphism with 
 respect to the group operation from point (a); 
 \item[(d)] the distance $d^{x}$ has the cone property with
respect to dilatations: foar any $u,v \in X$ such that $\displaystyle d(x,u)\leq 1$ and 
$\displaystyle d(x,v) \leq 1$  and all $\mu \in (0,A)$ we have: 
$$d^{x}(u,v) \ = \ \frac{1}{\mu} d^{x}(\delta_{\mu}^{x} u , \delta^{x}_{\mu} v) 
 \quad .$$ 
 \end{enumerate}
\label{tgene}
\end{thm}

The conical group $\displaystyle (U(x), \Sigma^{x}, \delta^{x})$ can be regarded as the tangent space 
of $(X,d, \delta)$ at $x$. 

By using proposition 5.4 \cite{siebert} and from some topological considerations
  we deduce the following characterisation of  tangent spaces asociated to some  
  dilatation structures. The following is corollary 4.7 \cite{buligacont}. 

\begin{cor}
Let $(X,d,\delta)$ be a dilatation structure with group $\Gamma = (0,+\infty)$
and the morphism $\nu$ equal to identity. 
Then for any $x \in X$ the local group  
 $\displaystyle (U(x), \Sigma^{x})$ is locally a simply connected Lie group 
 whose Lie algebra admits a positive graduation (a Carnot group).
 \label{cortang}
\end{cor}

\section{Linearity and dilatation structures}

In this section we describe the notion of linearity for dilatation structures, 
as in Buliga \cite{buligacont}.

\begin{defi}
Let $(X,d,\delta)$ be a    dilatation structure. A transformation $A:X\rightarrow X$ is linear if it 
is Lipschitz and it commutes with dilatations in the following sense: for any $x\in X$, $u \in U(x)$ and 
$\varepsilon \in \Gamma$, $\nu(\varepsilon) < 1$, if  $A(u) \in U(A(x))$ then  
$$ A \delta^{x}_{\varepsilon} = \delta^{A(x)} A(u) \quad  .$$
\label{defgl}
\end{defi}

In the particular case of  
$X$ finite dimensional real, normed vector space, 
$d$ the distance given by the norm, $\Gamma = (0,+\infty)$ and dilatations 
$\displaystyle \delta_{\varepsilon}^{x} u = x + \varepsilon(u-x)$, 
a linear transformations in the sense of definition \ref{defgl} is an affine transformation of the vector 
space $X$. More generally, linear transformations in the sense of definition 
\ref{defgl} have the expected 
properties related to linearity, as explained in section 5 \cite{buligacont}.


\begin{convention}
Further we shall say that  a property 
$\displaystyle \mathcal{P}(x_{1}$, $\displaystyle x_{2}$,  
$\displaystyle x_{3}, ...)$ holds for 
$\displaystyle x_{1}, x_{2}, x_{3},...$ sufficiently closed if for any 
compact, non empty set $K \subset X$, there is a positive constant $C(K)> 0$ 
such that $\displaystyle \mathcal{P}(x_{1},x_{2},x_{3}, ...)$ is true for any 
$\displaystyle x_{1},x_{2},x_{3}, ... \in K$ with 
$\displaystyle d(x_{i}, x_{j}) \leq C(K)$. 
\end{convention}

For example,  the expressions 
$$\delta_{\varepsilon}^{x} \delta^{y}_{\mu} z \quad , \quad
\delta^{\delta^{x}_{\varepsilon} y}_{\mu} \delta^{x}_{\varepsilon} z$$ 
are well defined for any $x,y,z \in X$ sufficiently closed and for any 
$\varepsilon,\mu \in 
\Gamma$ with $\nu(\varepsilon),\nu(\mu) \in (0,1]$. 
Indeed, let $K \subset X$ be compact, non empty set. Then there
is a constant $C(K) > 0$, depending on the set $K$ such that for any $\varepsilon,\mu \in 
\Gamma$ with $\nu(\varepsilon),\nu(\mu) \in (0,1]$ and any $x,y,z \in K$ with 
$d(x,y), d(x,z), d(y,z) \leq C(K)$  we have 
$$\delta^{y}_{\mu}z \in V_{\varepsilon}(x) \quad , \quad
\delta_{\varepsilon}^{x} z \in V_{\mu}(\delta^{x}_{\varepsilon} y) \quad .$$
Indeed, this is coming from the uniform (with respect to K) estimates:  
$$d(\delta^{x}_{\varepsilon} y, \delta^{x}_{\varepsilon} z) \leq \varepsilon 
d^{x}(y,z) + \varepsilon \mathcal{O}(\varepsilon) \quad , $$
$$d(x, \delta^{y}_{\mu} z) \leq d(x,y) + d(y, \delta^{y}_{\mu}z) 
\leq d(x,y) + \mu d^{y}(y,z) + \mu \mathcal{O}(\mu) \quad . $$
These estimates allow us to give the following definition.

\begin{defi}
A   dilatation structure $(X,d,\delta)$ is linear if for any $\varepsilon,\mu \in 
\Gamma$ such that $\nu(\varepsilon),\nu(\mu) \in (0,1]$, and for any 
$x,y,z \in X$ sufficiently closed we have 
$$\delta_{\varepsilon}^{x} \,  \delta^{y}_{\mu} z \ = \ 
\delta^{\delta^{x}_{\varepsilon} y}_{\mu} \delta^{x}_{\varepsilon} z \quad .$$
\label{defilin}
\end{defi}

Linear dilatation structures are very particular dilatation structures. 
The next theorem is theorem 5.7 \cite{buligacont}.  It is 
 shown that a linear and strong dilatation structure
comes from a normed conical group. 

\begin{thm}
Let $(X,d,\delta)$ be a linear   dilatation structure. Then the following 
two statements are equivalent: 
\begin{enumerate}
\item[(a)] For any $x \in X$ there is an open neighbourhood $D \subset X$
 of $x$ such that  $\displaystyle (\overline{D}, d^{x}, \delta)$ is the same 
 dilatation structure as the dilatation structure of the tangent space 
 of $(X,d,\delta)$ at $x$;   
\item[(b)]   The dilatation structure is strong (that is satisfies A4). 
\end{enumerate}
\label{tdilatlin}
\end{thm}

\section{Dilatation structures and inverse semigroups}

Here we prove that for dilatation structures linearity is equivalent 
to  a generalization of the statement from corollary \ref{corunu}. 
The result is new  for Carnot groups and the proof seems to be new even for 
vector spaces. 

\begin{defi}
A dilatation structure $(X,d,\delta)$ has the Menelaos property if 
for any two sufficiently closed $x,y \in X$  and for any 
$\varepsilon,\mu \in \Gamma$ with $\nu(\varepsilon),\nu(\mu) \in (0,1)$ we have 
$$\delta^{x}_{\varepsilon} \, \delta^{y}_{\mu} \ = \ \delta^{w}_{\varepsilon \mu}
\quad , $$
where $w \in X$ is the fixed point of the contraction $\displaystyle 
\delta^{x}_{\varepsilon} \delta^{y}_{\mu}$ (thus depending on $x,y$ and 
$\varepsilon , \mu$). 
\label{defmene}
\end{defi}

\begin{thm}
A linear dilatation structure has the Menelaos property. 
\label{thmenelin}
\end{thm}

\paragraph{Proof.} Let $x, y \in X$ be sufficiently closed and 
$\varepsilon,\mu \in \Gamma$ with $\nu(\varepsilon),\nu(\mu) \in (0,1)$. 
We shall define two sequences $\displaystyle x_{n}, y_{n} \in X$, $n \in
\mathbb{N}$. 

We begin with $\displaystyle x_{0} = x$, $y_{0} = y$. Let us define   by induction 
\begin{equation}
x_{n+1} \ = \ \delta_{\mu}^{\delta_{\varepsilon}^{x_{n}} y_{n}} x_{n} \quad , 
\quad y_{n+1} \ = \ \delta_{\varepsilon}^{x_{n}} y_{n} \quad .
\label{definduc}
\end{equation}

In order to check if the definition is correct we have to prove that for 
any $n \in \mathbb{N}$, if $\displaystyle x_{n}, y_{n}$ are
sufficiently closed then $\displaystyle x_{n+1}, y_{n+1}$ are
sufficiently closed too. 

Indeed, due to the linearity of the dilatation structure, 
 we can write the first part of (\ref{definduc}) as: 
$$x_{n+1} \ = \delta_{\varepsilon}^{x_{n}} \, \delta^{y_{n}}_{\mu} x_{n} \quad . $$
Then we can estimate the distance between $\displaystyle x_{n+1}, y_{n+1}$ like this: 
$$d(x_{n+1},y_{n+1}) \ = \ d( \delta_{\varepsilon}^{x_{n}} \, \delta^{y_{n}}_{\mu} x_{n}
, \delta_{\varepsilon}^{x_{n}} y_{n} ) \ = \ \nu(\varepsilon) \, 
d(\delta^{y_{n}}_{\mu} x_{n} , y_{n}) \ = \ \nu(\varepsilon \mu) d(x_{n}, y_{n})
\quad .$$
From $\nu(\varepsilon \mu) < 1$ it follows that $\displaystyle d(x_{n+1},
y_{n+1}) < d(x_{n}, y_{n})$, therefore  $\displaystyle x_{n+1}, y_{n+1}
$ are sufficiently closed.  We also find out that 
\begin{equation}
\lim_{n \rightarrow \infty} d(x_{n}, y_{n}) \ = 
\ 0 \quad . 
\label{equlim}
\end{equation}

Further  we use twice the linearity of the dilatation structure: 
$$\delta_{\varepsilon}^{x_{n}} \, \delta^{y_{n}}_{\mu} \ = \ 
\delta_{\mu}^{\delta_{\varepsilon}^{x_{n}} y_{n}} \, 
\delta_{\varepsilon}^{x_{n}} \ = \ \delta_{\varepsilon}^{\delta_{\mu}^{
\delta_{\varepsilon}^{x_{n}} y_{n}} x_{n}} \, 
\delta_{\mu}^{\delta_{\varepsilon}^{x_{n}} y_{n}} \quad . $$

By definition (\ref{definduc}) we arrive at the conclusion  that for any 
$n \in \mathbb{N}$ 
\begin{equation}
\delta_{\varepsilon}^{x_{n}} \, \delta^{y_{n}}_{\mu} \ = \
\delta_{\varepsilon}^{x} \, \delta^{y}_{\mu} \quad . 
\label{equlimain}
\end{equation}

From relation (\ref{equlimain}) we deduce that  the first part of 
(\ref{definduc}) can be written as: 
$$x_{n+1} \ = \delta_{\varepsilon}^{x_{n}} \, \delta^{y_{n}}_{\mu} x_{n} \ = \ 
\delta_{\varepsilon}^{x} \, \delta^{y}_{\mu} x_{n} \quad . $$
The transformation $\displaystyle \delta_{\varepsilon}^{x} \, \delta^{y}_{\mu}$ 
is a contraction of coefficient  $\nu(\varepsilon \mu) < 1$, therefore we easily
get: 
\begin{equation}
\lim_{n \rightarrow \infty} x_{n} \ = \ w \quad , 
\label{equlim2}
 \end{equation}
where $w$ is the unique fixed point of the contraction   $\displaystyle
\delta_{\varepsilon}^{x} \, \delta^{y}_{\mu}$. 

We put toghether (\ref{equlim}) and (\ref{equlim2}) and we get the limit: 
\begin{equation}
\lim_{n \rightarrow \infty} y_{n} \ = \ w \quad , 
\label{equlim3}
 \end{equation}
Using relations (\ref{equlim2}), (\ref{equlim3}), we may pass to the limit with $n \rightarrow \infty$ 
 in relation (\ref{equlimain}): 
 $$ \delta_{\varepsilon}^{x} \, \delta^{y}_{\mu} \ = \ \lim_{n \rightarrow 
 \infty} \delta_{\varepsilon}^{x_{n}} \, \delta^{y_{n}}_{\mu} \ = \ 
 \delta_{\varepsilon}^{w} \, \delta^{w}_{\mu} \ = \ \delta^{w}_{\varepsilon 
 \mu} \quad .$$ 
 The proof is done. \quad $\square$

\begin{cor}
Let $(X,d,\delta)$ be a strong linear dilatation structure, 
with group $\Gamma = (0,+\infty)$ and the morphism $\nu$ equal to identity. Any element of the  
inverse subsemigroup of $I(X)$ generated by  dilatations is locally 
a dilatation $\displaystyle \delta^{x}_{\varepsilon}$  or  a  left translation 
 $\displaystyle \Sigma^{x}(y, \cdot)$. 
\label{cordoi}
\end{cor}

\paragraph{Proof.}
Let $(X,d,\delta)$ be a strong linear dilatation structure. From the linearity 
and theorem \ref{thmenelin} we deduce that we have to care only about
the results of compositions of two dilatations $\displaystyle
\delta^{x}_{\varepsilon}$, $\displaystyle \delta^{y}_{\mu}$, with $\varepsilon
\mu = 1$.  
 
The dilatation structure is strong, therefore  by theorem \ref{tdilatlin} the
dilatation structure is locally coming from a conical group. In a conical group 
we can make the  following computation (here $\displaystyle \delta_{\varepsilon} =
\delta^{e}_{\varepsilon}$ with $e$ the neutral element of the conical group): 
$$ \delta^{x}_{\varepsilon} \delta^{y}_{\varepsilon^{-1}} z \ = \ 
x \delta_{\varepsilon} \left( x^{-1} y \delta_{\varepsilon^{-1}} \left(y^{-1} z 
\right)\right) \ = \ x \delta_{\varepsilon} \left( x^{-1} y \right) y^{-1} z
\quad . $$
Therefore the composition of dilatations $\displaystyle \delta^{x}_{\varepsilon}
\delta^{y}_{\mu}$, with $\varepsilon \mu = 1$, is a left translation. 

Another easy computation shows that composition of left translations with
dilatations are dilatations. The proof end by remarking that all the statements 
are local. \quad $\square$

\end{document}